\nonstopmode
\documentclass[a4paper, 10pt]{amsart}
\usepackage{latexsym}
\usepackage{fancyhdr}
\usepackage{amsmath, amssymb}
\usepackage[ansinew]{inputenc}
\usepackage[all]{xy}
\usepackage{verbatim}
\swapnumbers
\theoremstyle{plain}
\newtheorem*{lemma*}{Lemma}
\newtheorem{lemma}[subsection]{Lemma}
\newtheorem*{theorem*}{Theorem}
\newtheorem{theorem}[subsection]{Theorem}
\newtheorem*{proposition*}{Proposition}
\newtheorem{proposition}[subsection]{Proposition}
\newtheorem*{corollary*}{Corollary}
\newtheorem{corollary}[subsection]{Corollary}
\theoremstyle{definition}
\newtheorem*{definition*}{Definition}

\newtheorem*{example*}{Example}

\newtheorem*{algorithm*}{Algorithm}
\newtheorem*{remark*}{Remark}
\newtheorem{remark}[subsection]{Remark}

\newenvironment{demo}[1]{\par\smallskip\noindent{\bf #1.}}{\par\smallskip}
\numberwithin{equation}{subsection}

\def\ga{\gamma}
\def\de{\delta}
\def\ep{\epsilon}

\def\la{\lambda}

\def\rh{\rho}

\def\si{\sigma}

\def\ta{\tau}

\def\vh{\varphi}
\def\ch{\chi}

\def\De{\Delta}

\def\La{\Lambda}

\def\C{\mathbb{C}}

\def\I{\mathbb{I}}
\def\N{\mathbb{N}}

\def\R{\mathbb{R}}

\def\<{\langle}
\def\>{\rangle}
\renewcommand{\o}{\circ}

\let\on=\operatorname

\title[Perturbation of polynomials and operators]
{Perturbation of complex polynomials and normal operators}

\author[A. Rainer]
{Armin Rainer}

\address{Armin Rainer: Department of Mathematics, University of Toronto, 
40 St.\ George Street, Toronto, Ontario, Canada M5S 2E4}

\email{armin.rainer@univie.ac.at}

\begin{document}

\begin{abstract} 
We study the regularity of the roots of complex 
monic polynomials $P(t)$ of fixed degree depending smoothly on a real parameter $t$. 
We prove that each continuous parameterization of the roots of a generic $C^\infty$ curve $P(t)$ (which always exists)
is locally absolutely continuous. Generic means that no two of the continuously chosen roots meet of infinite order of flatness. 
Simple examples show that one cannot expect a better regularity than absolute continuity.
This result will follow from the proposition that for any $t_0$ there exists a positive integer $N$ such that 
$t \mapsto P(t_0\pm (t-t_0)^N)$ admits smooth parameterizations of its roots near $t_0$.
We show that $C^n$ curves $P(t)$ (where $n = \deg P$) admit differentiable roots 
if and only if the order of contact of the roots is $\ge 1$. 
We give applications to the perturbation theory of normal matrices and unbounded normal operators 
with compact resolvents and common domain of definition: 
The eigenvalues and eigenvectors of a generic $C^\infty$ curve of such operators can be arranged locally in 
an absolutely continuous way.
\end{abstract}

\thanks{The author was supported by the Austrian Science Fund (FWF), Grants P 17108-N04 \& J2771}
\date{January 21, 2009}

\maketitle

\section{Introduction}

Let us consider a curve of polynomials 
\[
P(t)(z) = z^n + \sum_{j=1}^n (-1)^j a_j(t) z^{n-j}, 
\]
where the coefficients $a_j : I \to \C$, $1 \le j \le n$, are complex valued functions defined on an interval $I \subseteq \R$. 
Given that the coefficients $a_j$ possess some regularity, it is natural to ask whether the roots 
of $P$ can be arranged in a regular way as well, i.e., whether it is possible to find $n$ regular functions 
$\la_j : I \to \C$, $1 \le j \le n$, such that $\la_1(t),\ldots,\la_n(t)$ represent the roots of $P(t)$ for each $t \in I$.

This problem has been extensively studied under the additional assumption that the polynomials $P(t)$ are hyperbolic, 
i.e., all roots of $P(t)$ are real. 
By a classical theorem due to Rellich \cite{Rellich37I}, there exist real analytic parameterizations of 
the roots of $P$ if its coefficients are real analytic. Bronshtein \cite{Bronshtein79} proved that     
if all $a_j$ are of class $C^n$, then there exists a differentiable parameterization of 
the roots of $P$ with locally bounded derivative (see also Wakabayashi \cite{Wakabayashi86} for a different proof).   
It has been shown in \cite{AKLM98} that if all $a_j$ are smooth ($C^\infty$) and no two of the increasingly ordered (thus continuous) 
roots meet of infinite order of flatness, then there exist smooth parameterizations of the roots. 
Moreover, by \cite{KLM04}, the roots may always be chosen twice differentiable provided that the $a_j$ are $C^{3n}$. 
The conclusion in this statement is best possible as shown by an example in \cite{BBCP06}. 
Recently, also the best possible assumptions have been found by \cite{ColombiniOrruPernazza08}: 
If the coefficients $a_j$ are $C^n$ (resp.\ $C^{2n}$), the roots allow $C^1$ (resp.\ twice differentiable) parameterizations. 
For further results on this problem see also 
\cite{Glaeser63R}, \cite{Dieudonne70}, \cite{Mandai85}, \cite{CC04}, \cite{LR07}, \cite{KurdykaPaunescu08}.

If the hyperbolicity assumption is dropped, then there is the following general result (e.g.\ \cite[II \S5 5.2]{Kato76}):
There exist continuous functions $\la_j : I \to \C$, $1 \le j \le n$, which parameterize the roots of a curve of polynomials $P$ 
with continuous coefficients $a_j : I \to \C$, $1 \le j \le n$. Note that, in contrast to the hyperbolic case, 
there is no hope that the roots of a polynomial $P$ which depends regularly on more than one 
parameter might be parameterized even continuously; just take $P(t,s)(z)=z^2-(t+i s)$, where $t,s \in \R$ and $i=\sqrt{-1}$. 
Of course, in that example the roots are given as $2$-valued analytic function with branching point $0$ 
in terms of a Puiseux series, e.g.\ \cite[Appendix]{Baumgaertel72}, but we do not go into that in this note. 
Here we restrict our attention to the one parameter case. 
In the absence of hyperbolicity the roots of a Lipschitz curve $t \mapsto P(t)$ of polynomials of degree $n$ 
may still be parameterized in a 
$C^{0,1/n}$ way, locally, which follows from a result of Ostrowski \cite{Ostrowski40},
but we cannot expect that the roots of $P$ are locally Lipschitz continuous 
even when the coefficients $a_j$ are real analytic; for instance, consider $P(t)(z) = z^2 - t$, $t \in \R$. However the roots of $P$ 
may possess a weaker regularity: They may be parameterized by locally absolutely continuous functions. In fact, Spagnolo 
\cite{Spagnolo99} proved that there exist absolutely continuous parameterizations of the roots of $P$ on compact intervals 
$I$ if one of the following conditions holds: 
\begin{enumerate}
\item $n=2$ and the coefficients $a_j$ belong to $C^5$,
\item $n=3$ and the coefficients $a_j$ belong to $C^{25}$ (the case $n=4$ is announced),
\item $P(t)(z) = z^n-f(t)$ and $f$ belongs to $C^{2n+1}$.
\end{enumerate}
The proof makes essential use of the explicit formulas for the roots of $P$ available in those cases.

In this paper we extend this result to generic smooth curves of polynomials $P$ of arbitrary degree $n$.
We say that $P$ is generic if no two of the continuously arranged roots of $P$ meet of 
infinite order of flatness.
We show in section \ref{AC} that, if the $a_j$ are smooth, then
there exists an absolutely continuous parameterization of the roots of $P$ on each compact interval $I$; 
actually, any continuous parameterization of the roots is locally absolutely continuous. 
In particular, these conditions are satisfied if the coefficients  $a_j$ are real analytic or, more generally,
belong to a quasianalytic class of $C^\infty$ functions. 
The main ingredient in the proof is the proposition \ref{4.3} that for any $t_0 \in I$ there exists a positive integer $N$ such that 
$t \mapsto P(t_0\pm (t-t_0)^N)$ admits smooth parameterizations of its roots near $t_0$.
It is not known whether the roots of $P$ may be arranged in a locally absolutely continuous way if $P$ is not generic.
That problem requires different methods. 

In section \ref{diff} we find conditions for the existence of differentiable parameterizations of 
the roots of $P$. Evidently, a necessary condition is that there exists a continuous choice of the roots 
such that whenever two of them meet they meet of order $\ge 1$.
We show that this condition is also sufficient, provided that the coefficients $a_j$ 
of $P$ belong to $C^n$.

In section \ref{reform} we discuss a reformulation of the problem in terms of a lifting problem which has been discussed in 
\cite{AKLM00} and \cite{KLMR05, KLMR06, KLMRadd}. This yields implicit sufficient conditions for a curve of 
polynomials $P$ to allow smooth, $C^1$, or twice differentiable parameterizations of its roots. 
As application we discuss the quadratic case.

Applications to the perturbation theory of normal matrices are given in section \ref{pert}. The eigenvalues and 
eigenvectors of a generic smooth curve $t \mapsto A(t)$ of normal complex matrices may be parameterized locally in an 
absolutely continuous way. The curve $t \mapsto A(t)$ is called generic if the associated characteristic polynomial 
$t \mapsto \ch_{A(t)}$ is generic. 
Examples show that without genericity or normality of $A(t)$ the eigenvectors need not admit continuous arrangements. 
We also prove that, for each $t_0$ there exists a positive integer $N$ 
such that $t \mapsto A(t_0\pm (t-t_0)^N)$ allows a smooth parameterization of its eigenvalues and eigenvectors near $t_0$. 
If $A$ is real analytic, then the eigenvalues and eigenvectors of $t \mapsto A(t_0\pm (t-t_0)^N)$ may be arranged real analytic as well.

In section \ref{pertuo} we obtain analogous results for curves $t \mapsto A(t)$ of unbounded 
normal operators in a Hilbert space with common domain of definition and with compact resolvents. 

For more on the perturbation theory of linear operators
consider Rellich \cite{Rellich37I, Rellich37II, Rellich39III, Rellich40IV, Rellich42V, Rellich69},
Kato \cite{Kato76}, Baumg\"artel \cite{Baumgaertel72}, and also \cite{AKLM98}, \cite{KM03}, and \cite{KMR}.

\nocite{RainerMONO}
\nocite{KLMR08}

\section{Preliminaries} \label{pre}

\subsection{\label{2.1}}

Let
\[
P(z) = z^n + \sum_{j=1}^n (-1)^j a_j z^{n-j} = \prod_{j=1}^n (z-\la_j)
\]
be a monic polynomial with coefficients $a_1,\ldots,a_n \in \mathbb{C}$ 
and roots $\la_1,\ldots,\la_n \in \mathbb{C}$. 
By Vieta's formulas, $a_i = \sigma_i(\la_1,\ldots,\la_n)$, where $\sigma_1,\ldots,\sigma_n$ 
are the elementary symmetric functions in $n$ variables:
\begin{equation} \label{esf}
\si_i(\la_1,\ldots,\la_n) = 
\sum_{1 \le j_1 < \cdots < j_i \le n} \la_{j_1} \cdots \la_{j_i}.
\end{equation}
Denote by $s_i$, $i \in \N$, the Newton polynomials 
$\sum_{j=1}^n \la_j^i$ which are related to the elementary symmetric functions 
by
\begin{equation} \label{rec}
s_k - s_{k-1} \sigma_1 + s_{k-2} \sigma_2 - \cdots 
+ (-1)^{k-1} s_1 \sigma_{k-1} + (-1)^k k \sigma_k = 0, \quad (k \ge 1).
\end{equation}

Let us consider the so-called Bezoutiant
\[ 
B := 
\begin{pmatrix} 
s_0 & s_1 & \ldots & s_{n-1}\\ 
s_1 & s_2 & \ldots & s_n \\ 
\vdots & \vdots & \ddots & \vdots\\ 
s_{n-1} & s_n & \ldots &  s_{2n-2} 
\end{pmatrix} 
= \left(s_{i+j-2}\right)_{1 \le i,j \le n}.
\]
Since the entries of $B$ are symmetric polynomials in $\la_1,\ldots,\la_n$, 
we find a unique symmetric $n \times n$ matrix $\tilde{B}$ with 
$B = \tilde{B} \circ \sigma$, where $\sigma =(\sigma_1,\ldots,\sigma_n)$. 

Let $B_k$ denote the minor formed by the first $k$ rows and columns of $B$. Then we find
\begin{equation} \label{eqdel}
\Delta_k(\la) := \on{det} B_k(\la) = \sum_{i_1 < i_2 < \cdots < i_k} 
(\la_{i_1}-\la_{i_2})^2 \cdots (\la_{i_1}-\la_{i_k})^2 \cdots (\la_{i_{k-1}}-\la_{i_k})^2.
\end{equation}
Since the polynomials $\Delta_k$ are symmetric, we have 
$\Delta_k = \tilde{\Delta}_k \circ \si$ 
for unique polynomials $\tilde{\Delta}_k$.

From \eqref{eqdel} follows that the number of
distinct roots of $P$ equals the
maximal $k$ such that $\tilde \De_k(P) \ne 0$.

\subsection{Multiplicity}  \label{2.2}

For a continuous real or complex valued function $f$ defined near $0$ in 
$\R$ let the \emph{multiplicity} (or \emph{order of flatness}) $m(f)$ at $0$ be the supremum of all integers 
$p$ such that $f(t) = t^p g(t)$ near $0$ for a continuous function $g$.
We define in the obvious way the \emph{multiplicity} $m_{t_0}(f)$ of $f$ at any $t_0 \in \R$ (if $f$ is defined near $t_0$).
Note that, if $f$ is of class $C^n$ and $m(f) < n$, then $f(t) = t^{m(f)} g(t)$ 
near $0$, where now $g$ is $C^{n-m(f)}$ and $g(0) \ne 0$. 

If $f$ is a continuous function on the space of polynomials, then for a fixed 
continuous curve $P$ of polynomials we will denote by $m(f)$ the multiplicity 
at $0$ of $t \mapsto f(P(t))$.

We shall say that two functions $f$ and $g$ \emph{meet of order $\ge p$ at $0$} 
or \emph{have order of contact $\ge p$ at $0$} iff $m(f-g) \ge p$.

\begin{lemma*} 
Let $I \subseteq \R$ be an interval containing $0$.
Consider a curve of polynomials
\[
P(t)(z) = z^n + \sum_{j=2}^n (-1)^j a_j(t) z^{n-j},
\]
with $a_j : I \to \C$, $2 \le j \le n$, smooth.
Then, for integers $r$, the following conditions are equivalent:
\begin{enumerate}
\item $m(a_k) \ge k r$, for all $2 \le k \le n$;
\item $m(\tilde{\Delta}_k) \ge k (k-1) r$, for all $2 \le k \le n$.
\end{enumerate}
\end{lemma*}

\begin{demo}{Proof}
Without loss of generality we can assume $r > 0$.

$(1) \Rightarrow (2)$: 
From \eqref{rec} we deduce $m(\tilde{s}_k) \ge k r$ for all $k$, where 
$s_k = \tilde{s}_k \circ \sigma$. By the definition of
$\tilde{\Delta}_k = \det(\tilde{B}_k)$ we obtain $(2)$.

$(2) \Rightarrow (1)$: 
It is easy to see that $\tilde \De_k(0)=0$ for all $2 \le k \le n$ implies $\tilde s_k(0) = 0$ for all $2 \le k \le n$.
Then by \eqref{rec} 
we have $a_k(0) = 0$ for all $2 \le k \le n$. 
So near $0$ we have $a_2(t) = t^{2 r} a_{2,2r}(t)$ and 
$a_k(t) = t^{m_k} a_{k,m_k}(t)$ for $3 \le k \le n$, 
where the $m_k$ are positive integers and $a_{2,2r}, a_{3,m_3},\ldots,a_{n,m_n}$ 
are smooth functions. 
Let us suppose for contradiction that for some $k > 2$ we have 
$m_k = m(a_k) < k r$. 
We put
\begin{equation} \label{m}
m := \min\left\{r,\frac{m_3}{3},\ldots,\frac{m_n}{n}\right\} < r,
\end{equation}
and consider the following continuous curve of polynomials for (small) 
$t \ge 0$:
\[
P_{(m)}(t)(z) := z^n + a_{2,2r}(t) t^{2 r - 2 m} z^{n-2}  
 - \cdots + (-1)^n a_{n,m_n}(t) t^{m_n-n m}.
\]
We have $\tilde \De_k(P_{(m)}(t)) = t^{-k(k-1)m} \tilde \De_k(P(t))$. 
By \eqref{m}, $r-m>0$, whence $t \mapsto \tilde \De_k(P_{(m)}(t))$, $2 \le k \le n$, vanish at $t=0$.
We may conclude as before that all 
coefficients of $t \mapsto P_{(m)}(t)$ vanish for $t=0$. 
But this is a contradiction for those $k$ with $m(a_k) = m_k = km$.
\qed\end{demo}

\begin{remark*}
If the coefficients $a_j$ of $P$ in lemma \ref{2.2} are just of class $C^n$, the conclusion remains 
true for $r = 1$. The proof is the same with the slight modification that  
we define $m_k := \min\{k,m(a_k)\}$ for all $k$. 
\end{remark*}

\subsection{Genericity condition} \label{2.3}

Let $I \subseteq \R$ be an interval. We call a curve of monic polynomials
\[
P(t)(z) = z^n + \sum_{j=1}^n (-1)^j a_j(t) z^{n-j}
\]
with continuous coefficients $a_j : I \to \C$, $1 \le j \le n$, \emph{generic} if the following  
equivalent conditions are satisfied at any $t_0 \in I$:
\begin{enumerate}
\item If two of the continuously parameterized roots of $P$ meet of infinite order of flatness at $t_0$, then
their germs at $t_0$ are equal.
\item Let $k$ be maximal with the property that the germ at $t_0$ of $t \mapsto \tilde \De_k(P(t))$ is not $0$. 
Then $t \mapsto \tilde \De_k(P(t))$ is not infinitely flat at $t_0$.
\end{enumerate}
The equivalence of $(1)$ and $(2)$ follows easily from \eqref{eqdel}.
For instance, $P$ is generic, if its coefficients are real analytic, or more generally, belong to a quasianalytic class of $C^\infty$ functions. 

\subsection{Splitting lemma} \cite[3.4]{AKLM98} \label{2.4}
{\it Let $P_0 = z^n + \sum_{j=1}^n (-1)^j a_j z^{n-j}$ be a polynomial satisfying 
$P_0 = P_1 \cdot P_2$, where $P_1$ and $P_2$ are polynomials without common root. 
Then for $P$ near $P_0$ we have $P = P_1(P) \cdot P_2(P)$ 
for real analytic mappings 
of monic polynomials $P \mapsto P_1(P)$ and $P \mapsto P_2(P)$, 
defined for $P$ 
near $P_0$, with the given initial values.}

\subsection{Absolutely continuous functions} \label{2.5}

Let $I \subseteq \R$ be an interval. A function $f : I \to \C$ is called 
\emph{absolutely continuous}, or $f \in AC(I)$, if for all $\ep > 0$ there exists a $\de >0$ 
such that $\sum_{i=1}^N (b_i-a_i) < \de$ implies $\sum_{i=1}^N |f(b_i)-f(a_i)| < \ep$,
for all sequences of pairwise disjoint subintervals $(a_i,b_i) \subseteq I$, $1 \le i \le N$.
By the fundamental theorem of calculus for the Lebesgue integral,
$f \in AC([a,b])$ if and only if there is a function $g \in L^1([a,b])$ such that 
\[
f(t) = f(a) + \int_a^t g(s) ds
\quad \text{for all } t \in [a,b]. 
\]
Then $f'=g$ almost everywhere.
Every Lipschitz function is absolutely continuous.

Gluing finitely many absolutely continuous functions provides an 
absolutely continuous function: Let $f_1 \in AC([a,b])$, $f_2 \in AC([b,c])$, and $f_1(b)=f_2(b)$.
Then $f : [a,c] \to \C$, defined by $f(t) = f_1(t)$ if $t \in [a,b]$ and $f(t)=f_2(t)$ if $t \in [b,c]$, belongs to $AC([a,c])$.
Similarly for more than two functions.

Let $\vh : I \to I$ be bijective, strictly increasing, and Lipschitz continuous. 
If $f \in AC(I)$ then also $f \circ \vh \in AC(I)$. Furthermore:

\begin{lemma*} 
Let $r >0$ and $n \in \N_{>0}$.
Let $f \in AC([0,r])$ (resp.\ $f \in AC([-r,0])$) and set $h(t) = f(\sqrt[n]{t})$
(resp.\ $h(t)=f(-\sqrt[n]{|t|})$). Then $h \in AC([0,r^n])$ (resp.\ $h \in AC([-r^n,0])$). 
\end{lemma*}

\begin{demo}{Proof}
There exists a function $g \in L^1([0,r])$ such that 
\[
f(t) = f(0) + \int_0^t g(s) ds
\]
for all $t \in [0,r]$. The function $(0,r^n] \to (0,r], t \mapsto \sqrt[n]{t}$, is smooth and bijective, so 
\[
\int_0^{r^n} |g(\sqrt[n]{s})| (\sqrt[n]{s})' ds = \int_0^r |g(s)| ds 
\]
and $t \mapsto g(\sqrt[n]{t}) (\sqrt[n]{t})'$ belongs to $L^1([0,r^n])$. Thus $h(t)=f(\sqrt[n]{t})$ is 
in $AC([0,r^n])$.

For the second statement consider the absolutely continuous function $f \o S|_{[0,r]}$, where $S : \R \to \R, t \mapsto -t$. 
By the above, $h_S(t) = (f \o S|_{[0,r]})(\sqrt[n]{t})$ is in $AC([0,r^n])$, and so 
$h(t) = h_S(S^{-1}|_{[-r^n,0]}(t)) = f(-\sqrt[n]{-t}) = f(-\sqrt[n]{|t|})$
is in $AC([-r^n,0])$. 
\qed\end{demo}

\section{Absolutely continuous parameterization of the roots} \label{AC}

\begin{lemma} \label{Lip}
Let $I \subseteq \R$ be an interval. Consider a curve of monic polynomials
\[
P(t)(z) = z^n + \sum_{j=1}^n (-1)^j a_j(t) z^{n-j}
\]
such that the coefficients $a_j : I \to \C$, $1 \le j \le n$, are continuous.
If there is a Lipschitz parameterization of the roots of $P(t)$, then any 
continuous parameterization is Lipschitz.
\end{lemma}

\begin{demo}{Proof}
Let $\mu_1,\ldots,\mu_n$ be a Lipschitz parameterization of the roots of $P$ on $I$ with common Lipschitz constant $C$.
Assume that $t \mapsto \la(t)$ is any continuous root of $t \mapsto P(t)$ for $t \in I$.
Let $t<s$ be in $I$.
Then there is an $i_0$ such that $\la(t)=\mu_{i_0}(t)$. Now let $t_1$ be the 
maximum of all $r\in [t,s]$ such that $\la(r)=\mu_{i_0}(r)$. If $t_1<s$ then 
$\mu_{i_0}(t_1)=\mu_{i_1}(t_1)$ for some $i_1\ne i_0$. 
Let $t_2$ be the 
maximum of all $r\in [t_1,s]$ such that $\la(r)=\mu_{i_1}(r)$. If $t_2<s$ then 
$\mu_{i_1}(t_2)=\mu_{i_2}(t_2)$ for some $i_2\notin \{i_0,i_1\}$. 
And so on until $s=t_k$ for some $k\le n$. Then we have (where $t_0=t$)
\begin{align*}
\frac{|\la(s)-\la(t)|}{s-t} 
&\le \sum_{j=0}^{k-1} \frac{|\mu_{i_j}(t_{j+1})-\mu_{i_j}(t_j)|}{t_{j+1}-t_j}
\cdot\frac{t_{j+1}-t_j}{s-t}
\le C.
\qed
\end{align*}
\end{demo}

\begin{proposition} \label{4.3}
Let $I \subseteq \R$ be an interval.
Consider a generic curve of monic polynomials
\[
P(t)(z) = z^n + \sum_{j=1}^n (-1)^j a_j(t) z^{n-j}, 
\]
with smooth coefficients $a_j : I \to \C$, $1 \le j \le n$. 
For any $t_0 \in I$, there exists a positive integer $N$ such that the roots of $t \mapsto P(t_0 \pm (t-t_0)^N)$ 
can be parameterized smoothly near $t_0$. If the coefficients $a_i$ are real analytic, then the roots of 
$t \mapsto P(t_0 \pm (t-t_0)^N)$ can be parameterized real analytically near $t_0$.
\end{proposition}

\begin{demo}{Proof}
It is no restriction to assume that $0 \in I$ and $t_0=0$.

We use the following:

\begin{algorithm*}
$(1)$ If all roots of $P(0)$ are pairwise different, the roots of $t \to P(\pm t)$ may be parameterized smoothly near $0$, 
by the implicit function theorem. Then $N=1$.

$(2)$ If there are distinct roots of $P(0)$, we put them into two subsets which factors $P(t) = P_1(t) P_2(t)$ by 
the splitting lemma \ref{2.4}. Suppose that $t \mapsto P_1(\pm t^{N_1})$ and $t \mapsto P_2(\pm t^{N_2})$ 
are smoothly solvable near $0$, then $t \mapsto P(\pm t^{N_1 N_2})$ is smoothly solvable near $0$ as well.

$(3)$ If all roots of $P(0)$ are equal, we reduce to the case $a_1=0$, by replacing $z$ with $z-a_1(t)/n$. Then all 
roots of $P(0)$ are equal to $0$, hence, $a_k(0) = 0$ for all $k$. 
Let $m:= \on{min}\{m(a_k)/k : 2 \le k \le n\}$ which exists since $P$ is generic (by lemma \ref{2.2}). 
Let $d$ be a minimal integer such that $d m \ge 1$.
Then for the multiplicity of $t \mapsto a_k(\pm t^d)$, $2 \le k \le n$, we find
\[
m(a_k(\pm t^d)) = d m(a_k) \ge d m k \ge k.
\] 
Hence we may write $a_k(\pm t^d) = t^k \tilde a_k^\pm(t)$ near $0$ with $\tilde a_k^\pm$ smooth, for all $k$.
Consider 
\[
\tilde P^\pm (t)(z) = z^n + \sum_{j=2}^n (-1)^j \tilde a_j^\pm(t) z^{n-j}.
\]
If $t \to \tilde P^\pm(t)$ is smoothly solvable and $t \mapsto \la_j^\pm(t)$ are its smooth roots, 
then $t \mapsto t \la_j^\pm(t)$ are smooth parameterizations of the roots of $t \mapsto P(\pm t^d)$.

Note that $m(\tilde a_k^\pm) = d m(a_k)-k$, for $2 \le k \le n$, and thus
\begin{equation} \label{tm}
\tilde m := \min_{2 \le k \le n}  \frac{m(\tilde a_k^\pm)}{k} = dm-1 < m, 
\end{equation}
by the minimality of $d$.

If $\tilde m = 0$ there exists some $k$ such that $\tilde a_k^\pm(0) \ne 0$, 
and not all roots of $\tilde P^\pm(0)$ are equal. 
We feed $\tilde P^\pm$ into step $(2)$.
Otherwise we feed $\tilde P^\pm$ into step $(3)$. 
\end{algorithm*}

Step $(1)$ and $(2)$ either provide a required parameterization or reduce the problem to a lower degree $n$.
Since $\tilde m$ is of the form $p/k$ where $2 \le k \le n$ and $p \in \N$ and by \eqref{tm}, 
also step $(3)$ is visited only finitely many times. 
So the algorithm stops after finitely many steps and it provides an integer $N$ and a smooth 
parameterization of the roots of $t \mapsto P(\pm t^N)$ near $0$. 
The real analytic case is analogous.
\qed\end{demo}

\begin{theorem} \label{4.4}
Let $I \subseteq  \R$ be an interval. Consider a generic curve of monic polynomials
\[
P(t)(z) = z^n + \sum_{j=1}^n (-1)^j a_j(t) z^{n-j}, 
\]
with smooth coefficients $a_j : I \to \C$, $1 \le j \le n$. 
Any continuous parameterization 
$\la=(\la_1,\ldots,\la_n) : I \to \C^n$
of the roots of $P$ is locally absolutely continuous. 
\end{theorem}

\begin{demo}{Proof}
It suffices to show that each $t_0 \in I$ has a neighborhood on which $\la$ is absolutely continuous.
Without restriction we assume that $0 \in I$ and $t_0=0$. 
By proposition \ref{4.3}, there is an integer $N$ and a neighborhood $J_N$ of $0$ such that 
$t \mapsto P(\pm t^N)$ allows a smooth parameterization $\mu^\pm=(\mu_1^\pm,\ldots,\mu_n^\pm)$ 
of its roots on $J_N$. 
Another continuous parameterization is provided by 
$t \mapsto \la(\pm t^N) = (\la_1(\pm t^N),\ldots,\la_n(\pm t^N))$.
By lemma \ref{Lip}, the parameterization $t \mapsto \la(\pm t^N)$ is actually Lipschitz (by shrinking $J_N$ if necessary), 
in particular, absolutely continuous.
Let $J = \{t \in I : \pm \sqrt[N]{|t|} \in J_N\}$, $J_{\ge 0} = \{t \in J : t \ge 0\}$, and $J_{\le 0} = \{t \in J : t \le 0\}$. 
By lemma \ref{2.5}, we find that $\la$ is absolutely continuous on $J_{\ge 0}$.
In order to see that $\la$ is absolutely continuous on $J_{\le 0}$ we apply lemma \ref{2.5} to $t \mapsto \la(- t^N)$,
if $N$ is even, and to $t \mapsto \la(t^N)$, if $N$ is odd.
Hence $\la$ is absolutely continuous on $J$.  
This completes the proof.
\qed\end{demo}

\begin{corollary} \label{4.5}
Any continuous parameterization of the roots of a real analytic, or more generally quasianalytic, 
curve $I \ni t \mapsto P(t)$ of monic polynomials is locally absolutely continuous. \qed 
\end{corollary}

\begin{remark}
The conclusion in theorem \ref{4.4} is best possible.
In general the roots cannot be chosen with first derivative in $L^p_{\on{loc}}$
for any $1 < p \le \infty$.
A counter example is given by
\[ 
P(t)(z) = z^n- t, \quad  t \in \R,
\]
if $n \ge \frac{p}{p-1}$, for $1 < p < \infty$,
and if $n \ge 2$, for $p=\infty$.

On the other hand, finding the optimal assumptions on $P$ for admitting locally absolutely continuous roots is an open problem.
\end{remark}

\section{Differentiable parameterization of the roots} \label{diff}

\begin{lemma} \label{4.1}
\cite[4.3]{KLMR05}
Consider a continuous curve $c : (a,b) \to X$ in a compact metric space $X$. 
Then the set of all accumulation points of $c(t)$ as 
$t \searrow a$ is connected.
\end{lemma}

\begin{proposition} \label{5.1}
Let $I \subseteq \R$ be an interval containing $0$.
Consider a curve of monic polynomials
\[
P(t)(z) = z^n + \sum_{j=1}^n (-1)^j a_j(t) z^{n-j}, 
\]
such that the coefficients $a_j : I \to \C$, $1 \le j \le n$, are of class $C^n$. 
Then the following conditions are equivalent:
\begin{enumerate}
\item There exists a local continuous
parameterization of the roots of $P$ 
near $0$ which is differentiable at $0$.
\item There exists a local continuous parameterization $\la_i$ of the roots of $P$ near $0$ such that 
$\la_i(0)=\la_j(0)$ implies $m(\la_i-\la_j) \ge 1$, for all $i \ne j$.
\item Split $P(t) = P_1(t) \cdots P_l(t)$ according to lemma \ref{2.4}, where $l$ is the number of 
distinct roots of $P(0)$. 
Then $m(\tilde \De_k(P_i)) \ge k(k-1)$, for all $1 \le i \le l$ and $2 \le k \le \deg P_i$.
\end{enumerate}
\end{proposition}

\begin{demo}{Proof} 
$(1) \Rightarrow (2)$ is obvious and $(2) \Rightarrow (3)$ follows immediately from \eqref{eqdel}.

$(3) \Rightarrow (1)$: Using the splitting $P(t) = P_1(t) \cdots P_l(t)$, we may suppose that all 
roots of $P(0)$ coincide. We can reduce to the case 
$a_1 = 0$ by replacing the variable $z$ with $z - a_1(t)/n$. 
Then all roots of $P(0)$ are equal to $0$.
By assumption and remark \ref{2.2}, we find $m(a_k) \ge k$, for all 
$2 \le k \le n$. So, for $t$ near $0$, we can write 
$a_k(t) = t^k a_{k,k}(t)$ for continuous $a_{k,k}$ and $2 \le k \le n$.
The continuous curve of polynomials 
$
P_{(1)}(t)(z) := 
z^n + \sum_{j=2}^n (-1)^j a_{j,j}(t) z^{n-j}
$ 
admits a continuous parameterization 
$\tilde \la = (\tilde \la_1,\ldots,\tilde \la_n)$ of its roots  
near $0$. Then $\la(t) := t \tilde \la(t)$ parameterizes the roots of $P$, locally near $0$, 
and is differentiable at $0$.
\qed\end{demo}

\begin{theorem} \label{5.2} 
Let $I \subseteq \R$ be an open interval.
Consider a curve of monic polynomials
\[
P(t)(z) = z^n + \sum_{j=1}^n (-1)^j a_j(t) z^{n-j}, 
\]
such that the coefficients $a_j : I \to \C$, $1 \le j \le n$, are of class $C^n$. 
Then the following conditions are equivalent:
\begin{enumerate}
\item There exists a global differentiable parameterization of the roots of $P$.
\item There exists a global continuous parameterization of the roots of $P$ with order of contact $\ge 1$ 
(i.e. if any two roots meet they meet of order $\ge 1$).
\item Let $t_0 \in I$. Split $P(t) = P_1(t) \cdots P_l(t)$ near $t_0$ according to lemma \ref{2.4}, where $l$ is the number of 
distinct roots of $P(t_0)$. 
Then $m_{t_0}(\tilde \De_k(P_i)) \ge k(k-1)$, for all $1 \le i \le l$ and $2 \le k \le \deg P_i$.
\end{enumerate}
\end{theorem}

\begin{demo}{Proof}
By proposition \ref{5.1}, it just remains to check $(3) \Rightarrow (1)$.

We use induction on $n$. 
There is nothing to prove if $n=1$.
So let us assume that $(3) \Rightarrow (1)$ holds for degrees strictly less 
than $n$. 

We may suppose that $a_1 = 0$ by replacing $z$ with $z - a_1(t)/n$. 
Consider the set $F$ of all $t \in I$ such that all roots of $P(t)$ coincide. 
Then $F$ is closed and its complement $I \backslash F$ is a
countable union of open subintervals whose boundary points lie in $F$. 

Let $J$ denote one such interval. For each $t_0 \in J$, the polynomial $P(t_0)$ has distinct 
roots which may be put into distinct subsets, and, by lemma \ref{2.4}, 
we obtain a local splitting $P(t) = P_1(t) P_2(t)$ near $t_0$, where 
both $P_1$ and $P_2$ have degree less than $n$. Clearly, $P_1$ and $P_2$ satisfy $(3)$ as well.
By induction hypothesis, we find differentiable parameterizations 
of the roots of $P$, locally near any $t_0 \in J$. 

Let $\la$ be a differentiable parameterization of the roots of $P$ defined on a maximal 
subinterval $J' \subseteq J$. Suppose that the right (say) endpoint $t_1$ of $J'$ 
belongs to $J$. Then there exists a differentiable parameterization 
$\bar \la$ of the roots of $P$, locally near $t_1$, and there is a 
$t_0 < t_1$ such that both $\la$ and $\bar \la$ are defined near $t_0$. 
Let $(t_m)$ be a sequence with $t_m \to t_0$. For each $m$ there exists 
a permutation $\ta_m$ such that $\la(t_m)=\ta_m.\bar \la(t_m)$. 
By passing to a subsequence, again denoted by $(t_m)$, 
we have $\la(t_m)=\ta.\bar \la(t_m)$ for a fixed permutation $\ta$ and for all $m$. 
Then $\la(t_0) = \lim_{t_m \to t_0} \la(t_m) = \ta.(\lim_{t_m \to t_0} \bar \la(t_m)) = \ta.\bar \la(t_0)$ and 
\[
\la'(t_0) = \lim_{t_m \to t_0} \frac{\la(t_m)-\la(t_0)}{t_m-t_0} 
= \lim_{t_m \to t_0} \frac{\ta.\bar \la(t_m)-\ta.\bar \la(t_0)}{t_m-t_0}
= \ta.\bar \la'(t_0).
\] 
Hence, the differentiable parameterization $\la$ of the roots of $P$ was not 
maximal: we can extend it differentiably by defining $\tilde \la(t) := \la(t)$ for $t \le t_0$ and 
$\tilde \la(t) := \ta.\bar \la(t)$ for $t \ge t_0$. This shows that there exists a differentiable parameterization 
$\la$ of the roots of $P$ defined on $J$.

Let us extend $\la$ to the closure of $J$, by setting it $0$ at the 
endpoints of $J$. Since $a_1 =0$, then $\la$ still parameterizes 
the roots of $P$ on the closure of $J$. 
Let $t_0$ denote the right (say) endpoint of $J$. 
By proposition \ref{5.1}, there exists a local continuous 
parameterization $\bar \la$ of the roots of $P$ near $t_0$ which is 
differentiable at $t_0$. 
Let $(t_m)$ be a sequence with $t_m \nearrow t_0$.  
By passing to a subsequence, we may assume that 
$\la(t_m)=\ta.\bar \la(t_m)$ for a fixed permutation $\ta$ and for all $m$. 
Then $\lim_{t_m \nearrow t_0} \la(t_m) 
= \ta.(\lim_{t_m \nearrow t_0} \bar \la(t_m)) = \ta.0 = 0$ and 
\[
\lim_{t_m \nearrow t_0} \frac{\la(t_m)}{t_m-t_0} 
= \lim_{t_m \nearrow t_0} \frac{\ta.\bar \la(t_m)}{t_m-t_0}
= \ta.\bar \la'(t_0).
\] 
It follows that the set of accumulation points of $\la(t)/(t-t_0)$, as 
$t \nearrow t_0$, lies in the $\on{S}_n$-orbit through $\bar \la'(t_0)$ of the symmetric group $\on{S}_n$. 
Since this orbit is finite, we may conclude from lemma \ref{4.1} 
that the limit $\lim_{t \nearrow t_0} \la(t)/(t-t_0)$ exits. 
Thus the one-sided derivative of $\la$ at $t_0$ exists. 

For isolated points in $F$, it follows from the discussion in the previous 
paragraph that we can apply a fixed permutation to one of the neighboring 
differentiable parameterizations of the roots in order to glue them 
differentiably. 
Therefore, we have found a differentiable parameterization $\la$ 
of the roots of 
$P$ defined on $I \backslash F'$, where $F'$ denotes the set of 
accumulation points of $F$.

Let us extend $\la$ by $0$ on $F'$. Then it provides a global differentiable 
parameterization of the roots of $P$, since any parameterization is 
differentiable at points $t' \in F'$. For: It is clear that the derivative 
at $t'$ of any differentiable parameterization has to be $0$. 
Let $\bar \la$ be the local parameterization near $t'$, provided by 
proposition \ref{5.1}. As above we may conclude that 
the set of accumulation points of $\la(t)/(t-t')$, as 
$t \to t'$, lies in the $\on{S}_n$-orbit through $\bar \la'(t') = 0$. 
\qed\end{demo}


\section{Reformulation of the problem} \label{reform}
\nocite{KLMRadd}

\subsection{Lifting curves over invariants} \label{6.1} 
Let $G$ be a compact Lie group and let $\rh : G \rightarrow \on{O}(V)$ be an
orthogonal representation in a real finite dimensional Euclidean
vector space $V$.
By a classical theorem of Hilbert and Nagata,
the algebra $\R[V]^{G}$ of invariant polynomials on $V$
is finitely generated.
So let $\si_1,\ldots,\si_n$ be a system of homogeneous generators
of $\R[V]^G$ of positive degrees $d_1,\ldots,d_n$.
Consider the orbit map
$\si = (\si_1,\ldots,\si_n) : V \to \R^n$.
The image $\si(V)$ is a semialgebraic subset of
$\{y \in \R^n : P(y) = 0 ~\mbox{for all}~ P \in I\}$,
where $I$ is the ideal of relations between $\si_1,\ldots,\si_n$.
Since $G$ is compact, $\si$ is proper and separates orbits of $G$,
it thus induces a homeomorphism between $V/G$ and $\si(V)$.

Let $H = G_v$ be the isotropy group of $v \in V$ and $(H)$ the conjugacy class 
of $H$ in $G$ which is called the type of the orbit $G.v$. The union $V_{(H)}$ 
of orbits of type $(H)$ is called an orbit type submanifold of the 
representation $\rh$, and $V_{(H)}/G$ is called an orbit type submanifold of 
the orbit space $V/G$. 
The collection of connected components of the manifolds $\{V_{(H)}/G\}$ 
forms a stratification of $V/G$ called orbit type stratification, 
see e.g.\ \cite[4.3]{Pflaum01A}.  

Let $c : \R \to V/G = \si(V) \subseteq \R^n$ be a smooth curve in the orbit space; smooth as curve in $\R^n$.
A curve $\bar c : \R \to V$ is called \emph{lift} of $c$ to $V$, 
if $c = \si \circ \bar c$ holds. 
The problem of lifting smooth curves over invariants is independent of the 
choice of a system of homogeneous generators of $\R[V]^G$, 
see \cite[2.2]{KLMR05}.  

Let $s \in \N$. 
Denote by $A_s$ the union of all strata $X$ of the orbit 
space $V/G$ with $\dim X \le s$, and by $I_s$ the ideal of 
$\R[V]^G$ consisting of all polynomials 
vanishing on $A_{s-1}$.
Let $c : \R \to V/G = \si(V) \subseteq \R^n$ be
a smooth curve, $t \in \R$, and $s = s(c,t)$ a minimal integer 
such that, for a neighborhood $J$ of $t$ in $\R$, 
we have $c(J) \subseteq A_s$. 
The curve $c$ is called \emph{normally nonflat at $t$} 
if there is $f \in I_s$ such that $f \circ c$ is not infinitely flat at $t$.
A smooth curve $c : \R \to \si(V) \subseteq \R^n$ is 
called \emph{generic}, if $c$ is normally nonflat at $t$ for each 
$t \in \R$.

Let $G = \on{S}_n$, the symmetric group, and let $\rh$ be the standard representation of 
$\on{S}_n$ in $\R^n$ by permuting the coordinates. The elementary symmetric functions $\si_i$ in \eqref{esf} 
generate the algebra of symmetric polynomials $\R[\R^n]^{\on{S}_n}$. Hence the image $\si(\R^n)$ 
may be identified with the space of monic hyperbolic polynomials of degree $n$.
Recall that a polynomial is called hyperbolic if all its roots are real. 
A lift to $\R^n$ of a curve $P$ in $\si(\R^n)$ represents a parameterization of the roots of $P$.
A curve $P$ of hyperbolic polynomials is generic in the sense of the last paragraph
if and only if it is generic in the sense of \ref{2.3}, see e.g.\ \cite[2.6]{LR07}. 

The following theorem generalizes the main results on the one dimensional perturbation theory of hyperbolic polynomials.
It collects the main results of \cite{AKLM00} and \cite{KLMR05, KLMR06, KLMRadd}.

\begin{theorem*}
Let $c : \R \to V/G = \si(V) \subseteq \R^n$ be a 
curve in the orbit space and let $d = \max\{d_1,\ldots,d_n\}$. Then:
\begin{enumerate}
\item If $c$ is real analytic, then it allows a real analytic lift, locally.
\item If $c$ is smooth and generic, then there exists a global smooth lift.
\item If $c$ is $C^d$, then there exists a global differentiable lift.
\end{enumerate}
If $G$ is finite, write $V = V_1 \oplus \cdots \oplus V_l$ as 
orthogonal direct sum of irreducible subspaces $V_i$ and define 
$k=\max\{d,k_1,\ldots,k_l\}$, where 
$k_i = \min\{|G.v| : v \in V_i \backslash \{0\}\}$. Then:
\begin{enumerate}
\item[$(4)$] If $c$ is $C^k$, then each differentiable lift is $C^1$.
\item[$(5)$] If $c$ is $C^{d+k}$, then there exists a global twice 
differentiable lift.
\end{enumerate}
\end{theorem*}

\subsection{\label{6.2}} 
Let us consider the standard action of the symmetric group $\on{S}_n$ on $\C^n$ 
by permuting the coordinates and the diagonal action of $\on{S}_n$ on 
$\R^n \times \R^n$ by permuting the coordinates in each factor.
Write $z = (z_1,\ldots,z_n) \in \C^n$ where $z_k = x_k + i y_k$, 
$1 \le k \le n$, $x=(x_1,\ldots,x_n) \in \R^n$, and 
$y=(y_1,\ldots,y_n) \in \R^n$. The mapping 
\[
T : \C^n \longrightarrow \R^n \times \R^n : 
z \longmapsto (x,y)
\]
is an equivariant $\R$-linear homeomorphism. Consequently, it descends 
to a homeomorphism $\hat T$ such that the following diagram commutes
\begin{equation} \label{diag}
\xymatrix{
\C^n \ar[rr]^{T} \ar[d] && \R^n \times \R^n \ar[d] \\
\C^n/\on{S}_n \ar[rr]_{\hat T} &&  (\R^n \times \R^n)/\on{S}_n 
}
\end{equation}
Consider the respective orbit type stratifications of the $\on{S}_n$-modules 
$\C^n$ and $\R^n \times \R^n$ and of its orbit spaces.
It is evident that $T$, and thus also $\hat T$, maps strata onto strata.
Note that, while the orbit type stratification of $\C^n/\on{S}_n \cong \C^n$ 
is finer than its stratification as affine 
variety, the orbit type stratification of $(\R^n \times \R^n)/\on{S}_n$ 
is its coarsest stratification, e.g.\ \cite[4.4.6]{Pflaum01A}.  
 
Let $P : \R \to \C^n/\on{S}_n = \C^n$ be a curve of monic polynomials of degree $n$. 
Then $\hat T \circ P$ is a curve in $(\R^n \times \R^n)/\on{S}_n \subseteq \R^N$. 
It follows that $P$ allows a regular lift to $\C^n$, i.e., a regular 
parameterization of its roots, if and only if $\hat T \circ P$ allows a 
lift of the same regularity to $\R^n \times \R^n$. 
Theorem \ref{6.1} provides sufficient conditions for $\hat T \circ P$ to be liftable 
regularly, and hence for $P$ to admit regular parameterizations of its roots.

As generators for the algebra $\C[\C^n]^{\on{S}_n}$ we may choose the Newton 
polynomials $s_i(z) = \sum_{j=1}^n z_j^i$, for $1 \le i \le n$.
By the first fundamental theorem of invariant theory for $\on{S}_n$ (e.g.\ \cite[3.4.1]{Smith95}), 
the algebra $\R[\R^n \times \R^n]^{\on{S}_n}$ is generated by the polarizations of 
the $s_i$:
\[
\ta_{i,j}(x,y) = \sum_{k=1}^n x_k^i y_k^j, 
\quad (i,j \in \N, 1 \le i+j \le n).
\]  
We may then identify the orbit projections 
\[\C^n \longrightarrow \C^n/\on{S}_n \quad \mbox{and} \quad 
\R^n \times \R^n \longrightarrow (\R^n \times \R^n)/\on{S}_n\]
with the mappings 
\[s=(s_i) : \C^n \longrightarrow s(\C^n) = \C^n \quad \mbox{and} \quad 
\ta = (\ta_{i,j}) : \R^n \times \R^n \longrightarrow \ta(\R^n \times \R^n) 
\subseteq \R^N,
\] 
respectively. 
Here $N = \binom{n+2}{n}-1 = \frac{1}{2} n(n+3)$. 
The image $\ta(\R^n \times \R^n)$ is a semialgebraic 
subset of $\R^N$. Since it is homeomorphic with $s(\C^n) = \C^n$, its dimension 
is $2n$. It follows that there are at least $\frac{1}{2} n(n-1)$ independent non-trivial polynomial 
relations between the $\ta_{i,j}$.  

The homeomorphism $\hat T$ from the diagram \eqref{diag} is then determined by:
\[
\hat T^{-1} : \R^N \supseteq \ta(\R^n \times \R^n) \longrightarrow \C^n : 
(\ta_{i,j}) \longmapsto \left(\sum_{k=0}^m \binom{m}{k} i^k \tau_{m-k,k}\right)
_{1 \le m \le n}.
\]

\subsection{The quadratic case} \label{6.3}
Without loss it suffices to consider $P(t)(z) = z^2-f(t)$ with $f : I \to \C$.
Let us consider the curve 
$\hat T \circ P$ in $(\R^2 \times \R^2)/\on{S}_2$ 
whose coordinates $\ta_{i,j}(P)$ have to satisfy:
\begin{gather*}
\ta_{1,0}(P) = \ta_{0,1}(P) = 0,~
\ta_{2,0}(P)-\ta_{0,2}(P) = 2 \on{Re}(f),~
\ta_{1,1}(P) = \on{Im}(f),\\
\ta_{2,0}(P) \ta_{0,2}(P) = \ta_{1,1}^2. 
\end{gather*}
It is easy to compute
\begin{equation} \label{TP}
\hat T \circ P = (0,0,|f|+\on{Re}(f),|f|-\on{Re}(f),\on{Im}(f)).
\end{equation}
In the following a square root of $f$ is any function $g$ satisfying $g^2=f$.
Applying \ref{6.1} and \ref{6.2}, we obtain:
{\it 
\begin{enumerate}
\item If $f$ is smooth and nowhere infinitely flat and $|f|$ is smooth, then there exist smooth square roots of $f$.
\item If $f$ and $|f|$ are of class $C^4$, then there exist twice differentiable square roots of $f$.
\end{enumerate}
}
\noindent
Theorem \ref{4.4} and theorem \ref{5.2} give:
{\it 
\begin{enumerate}
\item[$(3)$] If $f$ is smooth and nowhere infinitely flat, then any continuous choice of the square roots of $f$
is locally absolutely continuous.
\item[$(4)$] Assume that $f$ is $C^2$. Then there exist differentiable square roots of $f$
if and only if $f$ vanishes of order $\ge 2$ at all its zeros.
\end{enumerate}
}
Let us assume that $f$ is real valued. Then \eqref{TP} reduces to:
\[
(\hat T \circ P)(t) = \left\{
\begin{array}{c@{\quad \mbox{if}~}l}
(0,0,2 f(t),0,0) & f(t) \ge 0,\\
(0,0,0,-2 f(t),0) & f(t) \le 0. 
\end{array}
\right.
\] 
Suppose further that $f$ is $C^2$ and that $f(t_0) =0$ implies $f'(t_0)=f''(t_0)=0$. 
It follows that $\hat T \circ P$ is of class $C^2$.
By \ref{6.1} and \ref{6.2}, there exist $C^1$ parameterizations of the square root of $f$. So:
{\it 
\begin{enumerate}
\item[$(5)$] If $f$ is real valued, $C^2$, and $f(t_0) = 0$ implies $f'(t_0)=f''(t_0)=0$, then there exist $C^1$
square roots of $f$.
\end{enumerate}
}
\noindent
Combining $(3)$ and $(5)$ we obtain:
{\it 
\begin{enumerate}
\item[$(6)$] If $f$ is real valued and smooth, then each continuous choice of square roots of $f$ is locally absolutely continuous.
\end{enumerate}
}

\section{Regular diagonalization of normal matrices} \label{pert}

\subsection{\!\!\label{genericm}}
Let $I \subseteq \R$ be an interval.
A smooth curve of normal complex $n \times n$ matrices $I \ni t \mapsto A(t)=(A_{ij}(t))_{1 \le i,j \le n}$ is called \emph{generic},
if $I \ni t \mapsto \ch_{A(t)}$ is generic, where $\ch_{A(t)}(\la) = \on{det}(A(t)-\la \I)$ is the characteristic 
polynomial of $A(t)$.

\begin{theorem} \label{7.1}
Let $I \subseteq \R$ be an interval. Let $I \ni t \mapsto A(t)=(A_{ij}(t))_{1 \le i,j \le n}$ be a generic smooth curve of 
normal complex matrices acting on a complex vector space $V = \C^n$. Then:
\begin{enumerate} 
\item For each $t_0 \in I$ there exists an integer $N$ such that $t \mapsto A(t_0\pm(t-t_0)^N)$ allows a smooth 
parameterization of its eigenvalues and eigenvectors near $t_0$. If $A$ is real analytic, 
then the eigenvalues and eigenvectors of $t \mapsto A(t_0\pm(t-t_0)^N)$ may be parameterized real analytically near $t_0$.
\item There exist locally absolutely continuous parameterizations of the eigenvalues and the eigenvectors of $A$.
\end{enumerate} 
\end{theorem}

\begin{demo}{Proof}
We adapt the proof of \cite[7.6]{AKLM98}.

By theorem \ref{4.4}, the characteristic polynomial 
\begin{align} \label{ch}
\ch_{A(t)}(\la) = \on{det}(A(t) - \la \I) &= \sum_{j=0}^n (-1)^{n-j} \on{Trace}(\La^j A(t)) \la^{n-j} \\
&= (-1)^n \Big(\la^n + \sum_{j=1}^n (-1)^j a_j(t) \la^{n-j}\Big) \nonumber
\end{align}
admits a continuous, locally absolutely continuous parameterization $\la_1,\ldots,\la_n$ of its roots. 
This shows the first part of $(2)$.

Let us show $(1)$. Without loss we may assume that $t_0=0$. 
By proposition \ref{4.3}, there is an integer $N_0$ such that the eigenvalues of $t \mapsto A(\pm t^{N_0})$ 
can be parameterized by smooth functions $t \mapsto \mu_j^\pm(t)$ near $0$. 
Consider the following algorithm:

$(a)$ Not all eigenvalues of $A(0)$ agree. Let $\nu_1,\ldots,\nu_l$ denote the pairwise distinct eigenvalues of $A(0)$
with respective multiplicities $m_1,\ldots,m_l$. Assume without loss that 
\begin{align*}
\nu_1 &= \mu_1^\pm(0)=\cdots=\mu_{m_1}^\pm(0), \\
\nu_2 &= \mu_{m_1+1}^\pm(0)=\cdots=\mu_{m_1+m_2}^\pm(0),\\
\cdots &\cdots \cdots \cdots \cdots \cdots \cdots \cdots \cdots \cdots\\
\nu_l &= \mu_{n-m_l}^\pm(0)=\cdots=\mu_{n}^\pm(0).
\end{align*}
This defines a partition into subsets of smooth eigenvalues such that, for $t$ near $0$, they do not meet each other 
if they belong to different subsets.
For $1 \le j \le l$ consider
\begin{align*}
V_t^{(j),\pm} &:= \bigoplus_{\{i ~:~ \nu_j=\mu_i^\pm(0)\}} \on{ker}(A(\pm t^{N_0})-\mu_i^\pm(t)) \\
&= \on{ker} \big(\circ_{\{i ~:~ \nu_j=\mu_i^\pm(0)\}} (A(\pm t^{N_0})-\mu_i^\pm(t))\big).
\end{align*}
Note that the order of the compositions in the above expression is not relevant.
So $V_t^{(j),\pm}$ is the kernel of a smooth vector bundle homomorphism $B^\pm(t)$ of constant rank, 
and thus is a smooth vector subbundle of the trivial bundle $(-\ep,\ep) \times V \to (-\ep,\ep)$. 
This can be seen as follows: Choose a basis of $V$ such that $A(0)$ is diagonal. By the elimination procedure 
one can construct a basis for the kernel of $B^\pm(0)$. For $t$ near $0$, the elimination procedure (with the same choices) 
gives then a basis of the kernel of $B^\pm(t)$. The elements of this basis are then smooth in $t$ near $0$.

It follows that it suffices to find smooth eigenvectors in each subbundle $V^{(j),\pm}$ separately, expanded in the
constructed smooth frame field. But in this frame field the vector subbundle looks again like a constant vector space. 
So feed each of these parts ($t \to A(\pm t^{N_0})$ restricted to $V^{(j),\pm}$, as matrix with respect to the frame field) 
into step $(b)$ below.  

$(b)$ All eigenvalues of $A(0)$ coincide and are equal to $a_1(0)/n$, according to \eqref{ch}. 
Eigenvectors of $A(t)$ are also eigenvectors of $A(t) - (a_1(t)/n) \I$, thus 
we may replace $A(t)$ by $A(t) - (a_1(t)/n) \I$ and assume that the first coefficient of 
the characteristic polynomial \eqref{ch} vanishes identically. Then $A(0)=0$.

If $A(t)=0$ for $t$ near $0$ we choose the eigenvectors constant.

Otherwise write $A_{ij}(t) = t^m A_{ij}^{(m)}(t)$, where $m := \min\{m(A_{ij}) : 1 \le i,j \le n\}$ which exists 
by assumption. 
It follows from \eqref{ch} that the characteristic
polynomial of $A^{(m)}(t)$ is 
\[
\ch_{A^{(m)}(t)}(\la) = (-1)^n \Big(\la^n + \sum_{j=2}^n (-1)^j t^{-mj} a_j(t) \la^{n-j}\Big),
\]
Hence $m(a_k) \ge mk$ for all $k$.
By proposition \ref{4.3}, there exists an integer $N_1$ such that $t \mapsto \ch_{A^{(m)}(\pm t^{N_1})}$ 
admits smooth parameterizations of its roots (eigenvalues of $t \mapsto A^{(m)}(\pm t^{N_1})$) for $t$ near $0$. 
Eigenvectors of $A^{(m)}(\pm t^{N_1})$ are also eigenvectors of $A(\pm t^{N_1})$. 
There exist $1 \le i,j \le n$ such that $A_{i j}^{(m)}(0) \ne 0$ and thus
not all eigenvalues of $A^{(m)}(0)$ are equal. Feed $t \mapsto A^{(m)}(\pm t^{N_1})$ into $(a)$.

By assumption, this algorithm stops after finitely many steps and shows $(1)$. 
The real analytic case is analogous.  

Now we finish the proof of $(2)$. By $(1)$, we find an integer $N$ such that $t \mapsto A(\pm t^N)$ allows 
smooth parameterizations $t \mapsto \mu_j^\pm(t)$ and $t \mapsto v_j^\pm(t)$ of its eigenvalues and eigenvectors near $0$. 
In a similar way as in the proof of theorem \ref{4.4}, we can compose $t \mapsto \mu_j^\pm(t)$ and $t \mapsto v_j^\pm(t)$ 
with $t \mapsto \sqrt[N]{t}$ and $t \mapsto -\sqrt[N]{|t|}$ in order to obtain absolutely continuous parameterizations 
of the eigenvalues and eigenvectors of $A$ near $0$.
\qed\end{demo}

\begin{remark*}
The condition that $A(t)$ is normal cannot be omitted. Any choice of eigenvectors of the following real analytic curve 
$A$ of $2 \times 2$ matrices has a pole at $0$. Hence there does not exist an integer $N$ such that $t \mapsto A(\pm t^N)$ 
allows regular eigenvectors near $0$.
\[
A(t) =  
\begin{pmatrix}
0 & 1 \\
t & 0
\end{pmatrix}.
\]
The following smooth curve $A$ of symmetric real matrices allows smooth eigenvalues, but the eigenvectors cannot 
be chosen continuously. This example (due to \cite[\S 2]{Rellich37I}) shows that the assumption that $A$ is generic 
is essential in theorem \ref{7.1}.
\[
A(t) = e^{-\frac{1}{t^2}} 
\begin{pmatrix}
\cos \frac{2}{t} & \sin \frac{2}{t} \\
\sin \frac{2}{t} & - \cos \frac{2}{t}
\end{pmatrix}, \quad A(0) = 0.
\]
\end{remark*}

\section{Perturbation of unbounded normal operators} \label{pertuo}

\begin{theorem} \label{8.1}
Let $I \subseteq \R$ be an interval. Let $I \ni t \mapsto A(t)$ be a generic smooth curve of 
unbounded normal operators in a Hilbert space with common domain of definition and with compact resolvents. 
Let $t_0 \in I$ and let $z_0$ be an eigenvalue of $A(t_0)$. 
Let $n$ be the multiplicity of $z_0$.
Then:
\begin{enumerate} 
\item 
There exists an integer $N$ such that the $n$ eigenvalues of $t \mapsto A(t_0\pm(t-t_0)^N)$ passing through $z_0$ and the 
corresponding eigenvectors allow smooth parameterizations, locally near $t_0$.
If $A$ is real analytic, 
then the $n$ eigenvalues of $t \mapsto A(t_0\pm(t-t_0)^N)$ passing through $z_0$ and its eigenvectors may be arranged 
real analytically, locally near $t_0$.
\item There exist locally absolutely continuous parameterizations of the $n$ eigenvalues of $A$ passing through $z_0$ and its eigenvectors,
locally near $t_0$.
\end{enumerate} 
\end{theorem}

That $A(t)$ is a smooth (resp.\ real analytic) curve of unbounded operators means the following: 
There is a dense subspace $V$ of the Hilbert space $H$ such that $V$ is the domain of definition 
of each $A(t)$, and such that each $A(t)$ is closed and $A(t)^*A(t)=A(t)A(t)^*$, 
where the adjoint operator $A(t)^*$ is defined as usual by $\<A(t)u,v\>=\<u,A(t)^*v\>$ 
for all $v$ for which the left-hand side is bounded as function in $u \in H$. 
Note that the domain of definition of $A(t)^*$ is $V$.
Moreover, we require that $t \mapsto \<A(t)u,v\>$ is smooth (resp.\ real analytic) for each $u \in V$ and $v \in H$.
This implies that $t \mapsto A(t)u$ is of the same class $\R \to H$ for each $u \in V$, 
by \cite[2.3]{KM97} or \cite[2.6.2]{FK88}.

We call the curve $I \ni t \mapsto A(t)$ \emph{generic}, if no two unequal continuously parameterized 
eigenvalues meet of infinite order at any $t \in I$.

\begin{demo}{Proof}
We use the resolvent lemma in \cite{KM03} (see also \cite{AKLM98}):
{\it If $A(t)$ is smooth (resp.\ real analytic), then also the resolvent $(A(t)-z)^{-1}$ is smooth (resp.\ real analytic) into $L(H,H)$ 
in $t$ and $z$ jointly.}

Let $z$ be an eigenvalue of $A(s)$ of multiplicity $n$ for $s$ fixed. 
Choose a simple closed curve $\ga$ in the resolvent set of $A(s)$ enclosing only $z$ among all eigenvalues of $A(s)$.
Since the global resolvent set $\{(t,z) \in \R \times \C : (A(t)-z) : V \to H ~\text{is invertible}\}$ is open, 
no eigenvalue of $A(t)$ lies on $\ga$, for $t$ near $s$. Consider 
\[
t \mapsto -\frac{1}{2 \pi i} \int_\ga (A(t)-z)^{-1} dz =:P(t),  
\]
a smooth (resp.\ real analytic) curve of projections (on the direct sum of all eigenspaces corresponding to eigenvalues in the interior 
of $\ga$) with finite dimensional ranges and constant ranks (see \cite{AKLM98} or \cite{KM03}).
So for $t$ near $s$, there are equally many eigenvalues in the interior of $\ga$. Let us call them $\la_i(t)$, 
$1 \le i \le n$, (repeated with multiplicity) and let us denote by $e_i(t)$, $1 \le i \le n$, 
a corresponding system of eigenvectors of $A(t)$. 
Then by the residue theorem we have
\[
\sum_{i=1}^n \la_i(t)^p e_i(t) \< e_i(t), ~ \> = - \frac{1}{2 \pi i} \int_\ga z^p (A(t)-z)^{-1} dz 
\]
which is smooth (resp.\ real analytic) in $t$ near $s$, as a curve of operators in $L(H,H)$ of rank $n$.

Recall claim 2 in \cite[7.8]{AKLM98}:
{\it 
Let $t \mapsto T(t) \in L(H,H)$ be a smooth (resp.\ real analytic) curve of operators of rank $n$ in Hilbert space such that $T(0)T(0)(H)=T(0)(H)$. 
Then $t \mapsto \on{Trace}(T(t))$ is smooth (resp.\ real analytic) near $0$.
}

We conclude that the Newton polynomials 
\[
\sum_{i=1}^n \la_i(t)^p  = - \frac{1}{2 \pi i} \on{Trace}\int_\ga z^p (A(t)-z)^{-1} dz,
\]
are smooth (resp.\ real analytic) for $t$ near $s$,
and thus also the elementary symmetric functions
\[
\sum_{i_1 < \cdots < i_p} \la_{i_1}(t) \cdots \la_{i_p}(t).
\]
It follows that $\{\la_i(t) : 1 \le i \le n\}$ represents the set of roots of a polynomial of degree $n$ with smooth (resp.\ real analytic) 
coefficients. 
The statement of the theorem follows then from proposition \ref{4.3}, theorem \ref{4.4}, and theorem \ref{7.1}, since
the image of $t \mapsto P(t)$, for $t$ near $s$ describes a finite dimensional smooth (resp.\ real analytic) vector subbundle of 
$\R \times H \to \R$ and the $\la_i(t)$, $1 \le i \le n$, form the set of eigenvalues of 
$P(t)A(t)|_{P(t)(H)}$.
\qed\end{demo}

\def\cprime{$'$}
\providecommand{\bysame}{\leavevmode\hbox to3em{\hrulefill}\thinspace}
\providecommand{\MR}{\relax\ifhmode\unskip\space\fi MR }
\providecommand{\MRhref}[2]{%
  \href{http://www.ams.org/mathscinet-getitem?mr=#1}{#2}
}
\providecommand{\href}[2]{#2}


\end{document}